\documentclass[12pt]{amsart}
\usepackage{amsfonts, amsmath, amsthm,amscd,amssymb,verbatim,epsf}

\usepackage{psfig}
\usepackage{amsmath}
\usepackage{amscd}

\newtheorem{pro}{Proposition}[section]
\newtheorem{thm}[pro]{Theorem}
\newtheorem{lem}[pro]{Lemma}

\theoremstyle{definition}
\newtheorem{dfn}[pro]{Definition}

\theoremstyle{remark}

 \def\a{{\alpha}}

 \def\ra{{\rightarrow}}
 
 \def\g{{\gamma}}

 \def\2{{\mathbb Z_2}}

\def\a{\alpha}

\def\g{\gamma}

\def\wh{\widehat}
\def\wt{\widetilde}

\begin{document}

\title{Thick surfaces in hyperbolic 3-manifolds}
\author{Joseph D. Masters}

\begin{abstract}
We show that every closed, virtually fibered hyperbolic 3-manifold contains
 immersed, quasi-Fuchsian surfaces with convex cores of arbitrarily
 large thickness.
 \end{abstract}

\maketitle

\section{Introduction}
  
 This paper is concerned with the geometry of surfaces in 3-manifolds.
  Let us fix attention on a closed, hyperbolic 3-manifold,
 $M = \mathbb{H}^3/\Gamma$.  A famous, unproven conjecture of Waldhausen
 would imply that $\Gamma$ contains a surface subgroup. 
 By work of Marden, Thurston, and Bonahon, these
 come in two varieties:
 geometrically infinite, which are virtual fiber groups,
 and geometrically finite, which are quasi-Fuchsian.
 In the latter case, information about the geometric complexity
 of the group can have topological consequences.
 For example, $\pi_1$-injective immersions of surfaces can often
 be cut-and-pasted, or (in the cusped case) tubed around a boundary torus,
 to form new $\pi_1$-injective immersions; these procedures are often
 easier to carry out when the surface groups in question 
 are ``close'' to Fuchsian groups  (see \cite{CL}, \cite{M}, \cite{Wu}).

  Thus, assuming that $\Gamma$ contains surface groups,
 it is of interest to understand their geometric complexity.  The difficulty is
 that, given a topological description of an immersed surface, it is often
 difficult to describe the induced geometry.

 We shall show that if $M$ is fibered,
 and the surface is immersed transverse to the suspension flow,
 then the geometry of the immersion depends on the geometry
 of a certain polygon in $\mathbb{H}^2$,
 which was introduced by Cooper, Long and Reid.
 We are then able to
 show that if $M$ is virtually fibered,
 then $\Gamma$ contains quasi-Fuchsian surface subgroups
 of unbounded geometric complexity.
 To make a precise statement, we require some definitions.

  A 3-manifold is \textit{virtually fibered} if it has a finite cover
 which fibers over $S^1$.  We let $\mathbb{H}^n$ denote hyperbolic $n$-space,
 and we denote the hyperbolic metric by $d_{\mathbb{H}^n}$, or just $d$
 when the meaning is clear from context.  The boundary at infinity of
 $\mathbb{H}^n$ is denoted $S^{n-1}_{\infty}$.
 For the remainder of the paper  $M = \mathbb{H}^3/\Gamma$
 will represent a closed, virtually fibered, hyperbolic 3-manifold.
 In this setting, it will be convenient to define
 a \textit{surface group} to be the fundamental
 group of a \textbf{closed}, orientable surface of positive genus. 
 A Kleinian surface group $G$ is said to be \textit{quasi-Fuchsian}
 if its limit set is a topological circle.
 If $S$ is a closed, orientable surface of positive genus,
 we say that an immersion
 $f:S \rightarrow M$ is quasi-Fuchsian if it is $\pi_1$-injective,
 and $f_* \pi_1 S$ is a quasi-Fuchsian subgroup of $\Gamma$.

 Let $G \subset PSL_2(\mathbb{C})$ be a quasi-Fuchsian surface group,
 let $\Lambda(G)$ denote the limit set of $G$,
 let $Hull(\Lambda(G))$ denote the convex hull of $\Lambda(G)$,
 and let $Core(G) = Hull(\Lambda(G))/G$ be the convex core of $G$.
 If $G$ is not Fuchsian, then $\partial Core(G)$ has two components,
 $\partial^{\pm} Core(G)$, and the \textit{thickness} of $Core(G)$,
 denoted $t(G)$, is the distance between them in the hyperbolic metric;
 \textit{i.e}
 $$t(G)
 = \inf\{ d(x,y): x \in \partial^+ Core(G), y \in \partial^- Core(G)\}.$$
 It is also natural to define $t(G) = 0$ if $G$ is
 Fuchsian, and $t(G) = \infty$ if $G$ is geometrically infinite.

\begin{thm} \label{main}
If $M = \mathbb{H}^3/\Gamma$ is a closed hyperbolic 3-manifold which
 virtually fibers over the circle,
  then $\Gamma$ contains a sequence of quasi-Fuchsian
 surface subgroups $G_i$, such that $t(G_i) \rightarrow \infty$.
\end{thm}

\noindent
\textit{Remarks:}\\
1. It was proved in \cite{CLR} that every such $\Gamma$ contains
 a quasi-Fuchsian surface subgroup.\\
\\
2. In \cite{M}, examples were constructed of co-compact Kleinian groups
 $\Gamma$, for which the number of maximal, genus $g$, surface subgroups
 grows factorially with $g$.  However, the thicknesses of these subgroups
 remain bounded.\\
\\
3. We do not know if the converse is true.
\\
\\
\noindent
\textbf{Outline of proof}

We shall construct an explicit sequence
 of immersions into $M$, based on the ``cut and cross-join''
 construction of Cooper, Long and Reid.
 Lemma \ref{coarse} gives a simple criterion for showing
 that certain immersions of this type are quasi-Fuchsian;
 ours is essentially just a coarsening of a criterion
 given in \cite{CLR}.  Lemma \ref{thick} then
 gives a simple criterion, in terms of the ``leaf polygon'',
 for showing that a sequence of such immersions
 is thick, in the sense of Theorem \ref{main}.

 In the final section, we construct a sequence satisfying both
 criteria simultaneously.
 The main tool required here is the LERF property for surface groups.
 A detailed outline of the construction precedes the actual proof.

\section{General facts about hyperbolic geometry}
 \label{hypgeom}

For the proof of Theorem \ref{main}, we shall require some
 elementary facts about hyperbolic geometry and convex hulls-- see \cite{EM}
 for more background.  The first main goal is
 Lemma \ref{dense}, which characterizes thickness in terms of
 limit sets.

 Suppose we are given
 a sequence of quasi-Fuchsian surface groups $G_i$
 with limit sets $\Lambda_i$, convex hulls
 $Hull_i$, and convex cores $Core_i = Hull_i/G_i$.
  For each $i$ we may fix
 a homotopy equivalence $g_i:S_i \rightarrow Core_i$, where $S_i$ is a closed
 surface, and a lift $\tilde{g}_i:\tilde{S}_i \rightarrow \mathbb{H}^3$.  

 For each $x \in \mathbb{H}^n$,
 we let $S_x$ be the unit sphere
 in $T_x(\mathbb{H}^n)$, with standard round metric $d_{round}$,
 and we define a map 
 $\pi_x: S_x \rightarrow S^{n-1}_{\infty}$ by
 setting $\pi_x(v)$ to be the endpoint (at infinity) of the infinite 
 geodesic ray, based at $x$, with tangent vector $v$.
 Then, for any $p,q \in S^{n-1}_{\infty}$, we define
 $d_x(p,q) = d_{round}(\pi_x^{-1}p, \pi_x^{-1} q)$;
 the metric $d_x$ is called the \textit{visual metric}
 on $S^{n-1}_{\infty}$, based at $x$.
 If $(X, d)$ is any metric space, and $Y \subset X$,
 we let $N_r(Y, d)$ denote the $r$-neighborhood
 of $Y$, with respect to the metric $d$.

  We let $\overline{\mathbb{H}}^n= \mathbb{H}^n \cup S^{n-1}_{\infty}$.
 If $X \subset \mathbb{H}^n$ is any set, we let $\overline{X}$
 denote its closure in $\overline{\mathbb{H}}^n$, and let
 $X_{\infty} = \overline{X} \cap S^{n-1}_{\infty}$.  Note that, if
 $\tau$ is an isometry of $\mathbb{H}^n$, and $\overline{\tau}$
 is the extension of $\tau$ to $\overline{\mathbb{H}}^n$, then
 for any $x^{\prime} \in \mathbb{H}^n$ and
 $y^{\prime},z^{\prime} \in S^{n-1}_{\infty}$, we have
 $d_{{\tau} x^{\prime}} (\overline{\tau} y^{\prime},\overline{\tau} z^{\prime})
  = d_{x^{\prime}}(y^{\prime},z^{\prime})$. 
 If $x \in \mathbb{H}^3$ and $y \in S^2_{\infty}$,
 we let $\vec{xy}$ denote the geodesic ray
 from $x$ to $y$.
 The following fact will be useful:

\begin{lem}\label{theta}
If $x, x^{\prime} \in \mathbb{H}^3$, with $d(x,x^{\prime}) = r_0$,
 and if $y,z \in S^2_{\infty}$, then
 $d_{x^{\prime}}(y,z) \leq e^{r_0} d_x(y,z)$.
\end{lem}

\begin{proof}
 Let $T$ (resp. $T^{\prime}$) be a triangle with vertices $x,y$ and $z$,
 (resp. $x^{\prime}, y, z$), and let $\theta$ (resp. $\theta^{\prime}$)
 be the angle at $x$ (resp. $x^{\prime})$.  Thus $d_x(y,z) = \theta$,
 and $d_{x^{\prime}}(y,z) = \theta^{\prime}$.
 Let $P$ be the hyperbolic plane spanning $T$, and
 let $\pi: \mathbb{H}^3 \rightarrow P$ be the nearest-point projection.\\
\\
\textit{Claim:} $d(x, \pi x^{\prime}) \leq d(x, x^{\prime})$,
 and $d_{\pi x^{\prime}}(y,z) \geq d_{x^{\prime}}(y,z)$.\\
Proof of claim:  The triangle with vertices $x, x^{\prime}$
 and $\pi x^{\prime}$ has a right angle at $\pi x^{\prime}$.
 The hypotenuse is longer than either leg, and so
 $d(x, x^{\prime}) \geq d(x, \pi x^{\prime})$.

 Let $P_y$ and $P_z$ be geodesic planes perpendicular to $P$,
 and containing the rays $\vec{x^{\prime} y}$ and $\vec{x^{\prime}z}$,
 respectively.  Then $d_{x^{\prime}}(y,z)$ is equal to
 the angle made between $\vec{x^{\prime}y}$ and $\vec{x^{\prime} z}$,
 which is less than or equal to the angle between $P_y$ and $P_z$,
 which is equal to $d_{\pi x^{\prime}}(y,z)$. This proves the claim.\\

 By the claim, we may assume that $T$ and $T^{\prime}$ are co-planar.
 See Fig. 1.
  
\begin{figure}[!ht]
{\epsfxsize=3in \centerline{\epsfbox{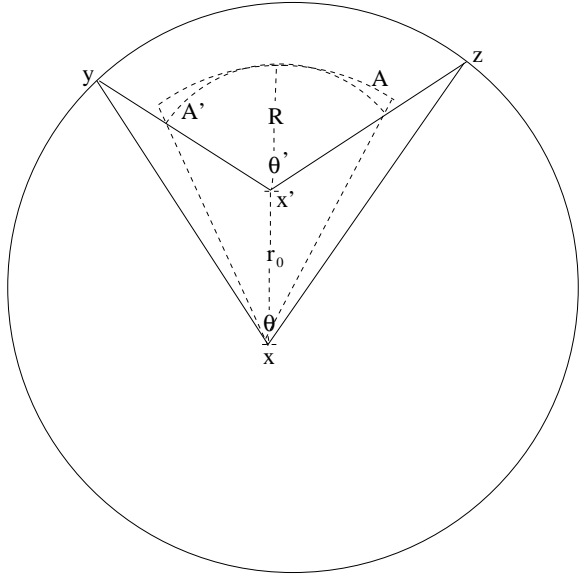}}\hspace{10mm}}
\caption{Difference between $d_x$ and $d_{x^{\prime}}$.}
\end{figure}

 Let $R>0$,
 let $B = N_{r_0+R}(x)$, and let $B^{\prime} = N_R(x^{\prime})$.
 Let $A^{\prime}$
 be a spherical geodesic on $\partial B^{\prime}$,
 connecting the two points of  $\partial B^{\prime} \cap T^{\prime}$,
 and let $A$ be the projection of $A^{\prime}$ onto $\partial B$.
 When written in polar coordinates centered
 at $x$, the hyperbolic metric in the plane of $T$ and $T^{\prime}$
 takes the form $(ds)^2 = (\sinh(r)  d \phi)^2
 + (dr)^2$, where $r$ represents the distance to $x$.
 We may parameterize the curve $A^{\prime}$ 
 as $(\phi(t), r(t))$, where $t \in [0,1]$,
 $(\phi(1/2), r(1/2)) = A \cap A^{\prime}$,
 $\phi(t)$ is monotone, and $r(t)$ is monotone
 on both $[0,1/2]$ and $[1/2, 1]$. 
 Note that $\Delta \phi \equiv |\phi(0) - \phi(1)| = \theta$, and that
 $\Delta r  \equiv Max(|r(0) - r(1/2)|, |r(1/2) - r(1)|) < 2R$.
 
One sees that
\begin{eqnarray*}
 Length(A^{\prime}) &<& \sinh(r_0+R) \Delta \phi + 2\Delta r\\ 
 &<& \sinh(r_0+R) \theta + 4R,
\end{eqnarray*}
 Also, using polar coordinates based at $x^{\prime}$,
 one sees that $Length(A^{\prime}) = \sinh(R) \theta^{\prime}$,
 and we have:

\begin{eqnarray*}
\theta^{\prime}
 &=& \frac{Length(A^{\prime})}{\sinh(R)}\\
 &<& \frac{\sinh(r_0+R) \theta + 4R}{\sinh(R)},
\end{eqnarray*}
 for every $R >0$.
It follows that
\begin{eqnarray*}
d_{x^{\prime}}(y,z) &=& \theta^{\prime}\\
 &\leq& \lim_{R \rightarrow \infty}\frac{\sinh(r_0+R) \theta + 4R}{\sinh(R)}\\
 &=& e^{r_0} \theta\\
 &=& e^{r_0} d_x(y,z).
\end{eqnarray*}
\end{proof}

 The following function will be used in the proof of
 Lemma \ref{dense}, and throughout the paper:\\

\begin{dfn} \label{rho} For any $\theta \in (0,\pi)$, we let
 $T(\theta)$ be the congruence class of a triangle in $\mathbb{H}^2$
 with two ideal vertices, and one non-ideal vertex
 of angle $\theta$.  We let 
 $\rho(\theta)$ be the distance from the non-ideal
 vertex of $T(\theta)$ to its opposite edge.
\end{dfn}
\vskip1pc
\begin{lem} \label{dense}
Suppose there is a sequence of numbers $\epsilon_i \rightarrow 0$,
 such that, for every $i$ and every $x \in \tilde{g}_i \tilde{S}_i$,
 we have $N_{\epsilon_i}(\Lambda_i, d_x) = S^2_{\infty}$.
  Then $t(G_i) \rightarrow \infty$.
\end{lem}

\begin{proof}
 We may assume that $\epsilon_i < \pi/2$, so $\Lambda_i$
 is not a round circle and $G_i$ is not Fuchsian.
 Since $g_i$ is a homotopy equivalence, then $g_i(S_i)$ separates
 the two boundary components of $Core_i$, and so $\tilde{g}_i \tilde{S}_i$
 separates the two boundary components of $Hull_i$. 
 Thus
\begin{eqnarray*}
t(G_i) \geq d(\tilde{g}_i\tilde{S_i}, \partial Hull_i).
\end{eqnarray*}

 Let $x \in \tilde{g}_i \tilde{S}_i$,
 and let $H$ be a half-space in $\mathbb{H}^3$
 such that $\overline{H}$ is disjoint from $\Lambda_i$.
 Since $N_{\epsilon_i}(\Lambda_i, d_x) = S^2_{\infty}$,
 and $H_{\infty} \cap \Lambda_i = \emptyset$,
 then the circle $\pi_x^{-1}H_{\infty}$
 has radius less than $\epsilon_i$ in the
 round metric on $S_x$.
 It then follows that the hyperbolic distance from $x$ to $H$ is at least
 $\rho(\epsilon_i)$, (see def \ref{rho}) and therefore
 \begin{eqnarray*} 
d(x, \partial Hull_i)
 &=& inf \{d(x, H): H \textrm{ is a half-space and }
 H_{\infty} \cap \Lambda_i = \emptyset \} \\
                      &\geq& \rho(\epsilon_i),
\end{eqnarray*} 
 for all $x \in \tilde{g}_i \tilde{S}_i$. Therefore
\begin{eqnarray*}
 t(G_i) &\geq& d(\tilde{g}_i \tilde{S}_i,\partial Hull_i) \\
       &\geq& \rho(\epsilon_i).
\end{eqnarray*}
 Then, since $\rho(x) \rightarrow \infty$ as $x \rightarrow 0$, we are done.
\end{proof}

We shall require the following technical variation on Lemma \ref{dense}.

\begin{lem} \label{technical}
Suppose there is a sequence of numbers $\epsilon_i \rightarrow 0$,
 a sequence of domains $D_i \subset \tilde{S}_i$,
 and a fixed number $R>0$, such that
 $N_R (G_i(\tilde{g}_i D_i), d_{\mathbb{H}^3}) \supset
 \tilde{g}_i \widetilde{S}_i$,
 and $N_{\epsilon_i}(\Lambda_i, d_x) = S^2_{\infty}$
 for every $x \in \tilde{g}_i D_i$.
  Then $t(G_i) \rightarrow \infty$.
\end{lem}

\begin{proof}
 Suppose we are given such a sequence, and
 let $x \in \tilde{g}_i \widetilde{S}_i$ and $y \in S^2_{\infty}$.
 Then there exists $\gamma \in G_i$ and $x^{\prime} \in \tilde{g}_i D_i$
 such that
 $d(\gamma x, x^{\prime}) \leq R$.
 
 We have:
\begin{eqnarray*}
d_x(\Lambda_i, y) &=&
  d_{\gamma  x}(\overline{\gamma} \Lambda_i,
                                  \overline{\gamma} y)\\
 &=& d_{\gamma  x}(\Lambda_i, \overline{\gamma} y)\\
 &\leq& e^R  d_{x^{\prime}}(\Lambda_i,
      \overline{\gamma} y) \,\,\, \textrm{  by Lemma \ref{theta}},\\
&\leq& e^R \epsilon_i.
\end{eqnarray*}
 So, letting $\epsilon_i^{\prime} = e^R \epsilon_i$,
 we obtain a sequence satisfying the conditions of Lemma \ref{dense},
 and so $t(G_i) \rightarrow \infty$.
\end{proof}

 We shall also require the existence of the Cannon-Thurston map.
 Let $M$ be a hyperbolic
 3-manifold fibering over $S^1$, with fiber $F$,
 let $\mathcal{I}:F \rightarrow M$ be the inclusion map,
 and let $\tilde{\mathcal{I}}: \mathbb{H}^2 \rightarrow \mathbb{H}^3$
 be a lift of $\mathcal{I}$.

\begin{thm} (Cannon-Thurston)
 The map  $\tilde{\mathcal{I}}$ extends continuously to a map
 $\overline{\mathcal{I}} :\overline{\mathbb{H}}^2 \rightarrow
 \overline{\mathbb{H}}^3$.
\end{thm} 

 This was proved in \cite{CT}.
 A published proof may be found in the appendix of \cite{CLR}.

 If $\gamma \in \pi_1 F$, then there are corresponding actions
 $\tau_{\gamma}$, of $\gamma$ on $\mathbb{H}^2$, and
 $\mathcal{I}_* \gamma$, of $\gamma$ on $\mathbb{H}^3$.
 One may check that the following diagram commutes:

\begin{center}
$\begin{CD}\overline{\mathbb{H}}^2 @>\overline{\tau}_{\gamma}>>  \overline{ \mathbb{H}}^2\\
@VV{\overline{\mathcal{I}}}V                @VV{\overline{\mathcal{I}}}V\\
\overline{\mathbb{H}}^3  @>\overline{\mathcal{I}_*\gamma}>>  \overline{\mathbb{H}}^3
\end{CD}$
\end{center}

\section{Immersions in bundles: Quasi-Fuchsian criterion}\label{immersions}
In \cite{CLR},
 Cooper, Long and Reid constructed a class of essential immersions,
 called ``cut-and-cross-join surfaces''
 into arbitrary hyperbolic surface bundles. 
 They also gave checkable criteria
 for determining whether such an immersion is geometrically finite or
 geometrically infinite.
 In this section, we review their construction and their criterion
 for geometric finiteness, and then
 give a (presumably) coarser criterion, which will be
 easier to verify for our examples.
 
 Suppose $M = \mathbb{H}^3/\Gamma$
 is a hyperbolic 3-manifold,
 which fibers over the circle
 with fiber $F$, and pseudo-Anosov monodromy $f:F \rightarrow F$.
 Representing $M$ as the mapping torus
 $F \times [0,1]/(x,1) = (fx, 0)$, the
 \textit{suspension flow} on $M$ has flow lines of the form
 $\bigcup_{i \in \mathbb{Z}}(f^i p \times [0,1])$.

 Let $\alpha \subset F$ be a simple, non-separating closed curve.
 Suppose $\pi:\widehat{F} \rightarrow F$ is a finite cover of $F$,
 containing 1-1 lifts $\widehat{\alpha}$ and $\widehat{f \alpha}$ of
 $\alpha$ and $f \alpha$, respectively, such that
 $\widehat{\alpha} \cap \widehat{f \alpha} = \emptyset$.
 Let $ \mathring{N}(\widehat{\alpha} \cup \widehat{f \alpha})$
 denote the interior of a regular neighborhood of
 $\widehat{\alpha} \cup \widehat{f \alpha}$, and let
 $\widehat{F}_0 = \widehat{F} - \mathring{N}(\widehat{\alpha} \cup \widehat{f \alpha})$.
 We label the components of
$\partial \widehat{F}_0$ as
 $\widehat{\alpha}^+,\widehat{\alpha}^-,  \widehat{f \alpha}^+$,
 and $ \widehat{f \alpha}^-$,
 where $f \pi(\widehat{\alpha}^+) = \pi(\widehat{f \alpha}^+)$.
 Let $g: \widehat{F}_0 \rightarrow F$ be the restriction of $\pi$;
 let $A^{\pm}$ be a pair of annuli, and let
 $g:A^{\pm} \rightarrow M$ be a map transverse to the suspension
 flow, such that
 $g|\partial A^{\pm}$ is a homeomorphism onto
 $\widehat{\alpha}^{\pm} \cup \widehat{f \alpha}^{\mp}$.
 We glue the annuli to $\widehat{F}_0$ by identifying
 $\partial A^{\pm}$
 with $\widehat{\alpha}^{\pm} \cup \widehat{f \alpha}^{\mp}$,
 to get a surface $S$ with the same genus as $\widehat{F}$.
 The maps $g$ which we have defined piece together
 to give a map $g: S \rightarrow M$ 
 which is transverse to the suspension flow, and therefore
 $\pi_1$-injective, by \cite{CLR} or \cite{M} or \cite{F}.
 Note that $S$ is connected if and only if
 $\widehat{\alpha} \cup \widehat{f \alpha}$ is non-separating.
 In the event that $S$ is disconnected, we replace $S$ with
 one of its connected components.

 We will use the notation $S = S(\widehat{\alpha}, \widehat{f\alpha})$
 for a surface $S$ constructed in this manner;
 the corresponding immersion is denoted
 $g:S(\widehat{\alpha}, \widehat{f\alpha}) \rightarrow M$.

 There is a useful criterion for showing that immersions of this
 type are quasi-Fuchsian.  We represent
 the universal cover of $M$ as
 $\widetilde{M} = \mathbb{H}^2 \times \mathbb{R}$,
 where $\mathbb{H}^2$ is the universal cover of $F$.
 Let $\widetilde{S}$ be the universal cover of $S$,
 and let $\tilde{g}: \widetilde{S} \rightarrow \widetilde{M}$
 be a lift of $g$.
 Let $\phi: \widetilde{M} \rightarrow \mathbb{H}^2$ be the map
 which identifies flow lines to points.
 Define $A(S)$ to be the closure of $\phi(\tilde{g}(\widetilde{S}))$ in
 $\mathbb{H}^2$.
  Then $A(S)$ is a convex subset of $\mathbb{H}^2$
 (\cite{CLR} Proposition 3.9).   We also have (\cite{CLR} Thm 3.14):

 \begin{thm} \label{gf} (Cooper-Long-Reid)
The immersion $g:S \rightarrow M$ is quasi-Fuchsian if and only
 if $A(S) \neq \mathbb{H}^2$.
 \end{thm}
 
 It will be useful to have a more explicit picture of $A(S)$.
 To that end, we let
 $S_0 = S - (\widehat{\alpha}^+ \cup \widehat{f \alpha}^+)$, let
 $\widetilde{S}_0 \subset \widetilde{S}$
 be a connected component
 of the pre-image of $S_0$, and let
 $A_0 = \overline{\phi(\tilde{g} (\widetilde{S}_0))}$. 
 Assume that $A_i$ has been constructed, and that each
 of the (countably many) components of $Frontier(A_i)$
 is a lift of $f^n \alpha$ for some $n$.
 Let $\{ \beta_{ij}: j = 1, 2, ... \}$ be the set of boundary
 components of $A_i$.
  For each $\beta_{ij}$, there is a map
 $\tau_{ij}:\mathbb{H}^2 \rightarrow \mathbb{H}^2$,
 lifting a power of $f$, such that
 $\tau_{ij} A_i \cap A_i = \beta_{ij}$, which is well-defined
 up to a covering translation which leaves $\tau_{ij} A_i$ invariant.
 Let $A_{i+1} = A_i \cup \,\,\, \bigcup_{j=1,2,...} \tau_{ij}A_i$.
 We then have $A_i$ defined recursively for all $i \geq 0$.
 Let $A(S)$ be the closure of $\bigcup_i A_i$.
 It can be checked that this agrees with the previous definition.

\begin{lem} \label{coarse}
 (quasi-Fuchsian criterion)
Suppose that some $\tau_{ij}$ has a fixed point in $\mathbb{H}^2$.
 Then $S$ is quasi-Fuchsian.
\end{lem}

\begin{proof}
 Let $\tilde{p} \in \mathbb{H}^2$ be a fixed point for $\tau_{ij}$
 (see Figure 2).
 We have, by assumption, $\tau_{ij} A_i \cap A_i = \beta_{ij}$.
 Since no simple closed curve in $F$ is invariant under $f$,
 then $\tau_{ij} \beta_{ij} \neq \beta_{ij}$.
 Thus $\tau_{ij} \beta_{ij} \cap \beta_{ij} = \emptyset$, and
 it follows that $\tilde{p} \not\in A_i$.

\begin{figure}[!ht]
{\epsfxsize=4in \centerline{\epsfbox{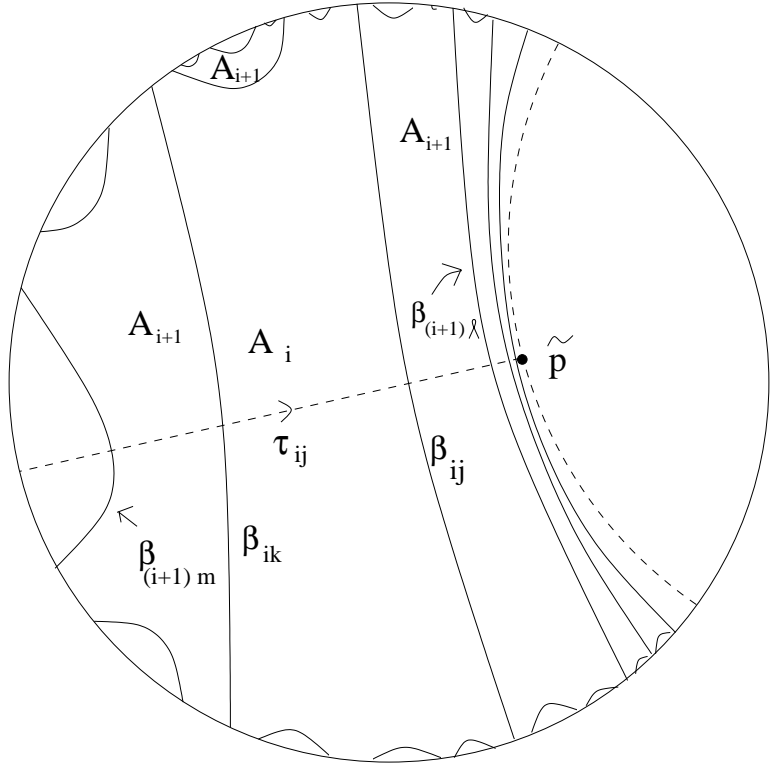}}\hspace{10mm}}
\caption{The existence of a fixed point
 implies that the surface group is quasi-Fuchsian.}
\end{figure}

 Note that $\tau_{ij}^{-1} A_i \cap A_i = \tau_{ij}^{-1} \beta_{ij}$,
 is a component $\beta_{i k}$ of $Frontier \, A_i$, and so
 we may take $\tau_{ij}^{-1} = \tau_{i k}$.
 Then $\beta_{(i+1) \ell} \equiv \tau_{ij} \beta_{ij}$
 and $\beta_{(i+1) m} \equiv \tau_{ik} \beta_{ik}$
 are both components of $Frontier A_{i+1}$.
 Since  $\tau_{ij}^3 \beta_{(i+1) m} = \tau_{ij}^3(\tau_{ik} \beta_{ik})
              = \tau_{ij}^3(\tau_{ij}^{-2} \beta_{ij})
              = \tau_{ij}(\beta_{ij})              
              = \beta_{(i+1) \ell}$,
 we may take $\tau_{ij}^3 = \tau_{(i+1) \ell}$, and 
 then $\tilde{p}$ is a fixed point for $\tau_{(i+1) \ell}$.
 Repeating our prior argument, we conclude that
 $\tilde{p} \not \in A_{i+1}$.  An inductive argument
 shows that $\tilde{p} \not \in A_j$ for all $j$, and
 so $int A(S) \neq \mathbb{H}^2$. Since $int A(S)$ is convex,
 we conclude that $A(S) \neq \mathbb{H}^2$.
 Therefore, by Theorem \ref{gf},
 the immersion $g:S \rightarrow M$ is quasi-Fuchsian.
\end{proof}

\section{Immersions in bundles: thickness criterion}

\subsection{Visual distortion of the Cannon-Thurston map} ${}$

Our first goal is to control the visual distortion of the Cannon-Thurston
 map, which was introduced in Section \ref{hypgeom}.
 We fix a point $x \in \mathbb{H}^2$, and consider the restriction
 $\overline{\mathcal{I}}:(S^1_{\infty}, d_x)
 \rightarrow (S^2_{\infty}, d_{\mathcal{I} x})$
 as a map of metric spaces.
 By compactness, $\overline{\mathcal{I}}$ is uniformly continuous,
 so there is a function $\sigma_x:(0,1] \rightarrow (0,1]$
 such that, given $\delta \in (0,1]$, and
 $y,z \in S^1_{\infty}$ satisfying $d_x(y,z) < \sigma_x(\delta)$,  
 then
 $d_{\mathcal{I} x}(\overline{\mathcal{I}} y, \overline{\mathcal{I}} z)
 < \delta$.
 Since the metrics $d_x$ and $d_{\mathcal{I} x}$ vary continuously
 with $x$, we may assume that
 the functions $\sigma_x$ vary continuously with $x$.

\begin{lem} \label{invarianceh}
 ($\pi_1 F$-invariance)
Let $x \in \mathbb{H}^2$. Then, given any $\gamma \in \pi_1 F$,
 any $\delta \in (0,1]$, and any $y,z \in S^1_{\infty}$
 such that $d_{\tau_{\gamma} x}(y,z) < \sigma_x(\delta)$, we have
 $d_{\mathcal{I} \, \tau_{\gamma} \, x}
 (\overline{\mathcal{I}} y, \overline{\mathcal{I}} z) < \delta$.
\end{lem}

\begin{proof}
 
 Given $\delta \in (0,1]$,
 $\gamma \in \pi_1 F$, and $y,z \in S^1_{\infty}$ such that
 $d_{\tau_{\gamma} x}(y,z) < \sigma_x(\delta)$. We have:

\begin{eqnarray*}
d_{\mathcal{I} \tau_{\gamma} x}
(\overline{\mathcal{I}} y, \overline{\mathcal{I}}z)
 &=& d_{\mathcal{I}_* \gamma^{-1} \,\, \mathcal{I} \tau_{\gamma} x}
(\overline{\mathcal{I}_* \gamma}^{-1}\overline{\mathcal{I}} y,
 \overline{\mathcal{I}_* \gamma}^{-1}\overline{\mathcal{I}} z)\\
 &=& d_{\mathcal{I} x} (\overline{\mathcal{I}}
 \overline{\tau}_{\gamma}^{-1} y,
 \overline{\mathcal{I}} \overline{\tau}_{\gamma}^{-1}z),
 \textrm{ by the commutative diagram from Section 2}.
\end{eqnarray*}

 We also have $d_x(\overline{\tau}_{\gamma}^{-1} y,
            \overline{\tau}_{\gamma}^{-1} z)= d_{\tau_{\gamma} x}(y,z)
         < \sigma_x(\delta)$, and so
  $d_{\mathcal{I}x}( \overline{\mathcal{I}} \overline{\tau}_{\gamma}^{-1} y,
           \overline{\mathcal{I}} \overline{\tau}_{\gamma}^{-1}z)
 < \delta$.
This proves the lemma.
\end{proof}
 Since $F$ is compact, then we may fix a compact fundamental
 domain $D$ for $F$ in $\mathbb{H}^2$, and
 by Lemma \ref{invarianceh}, we may assume that,
 for all $x \in D$ and all $\gamma \in \pi_1 F$, we have
 $\sigma_{\tau_{\gamma} x} = \sigma_x$.
 Then we define a function $\sigma:(0,1] \rightarrow (0,1]$ by the rule:
 $\sigma(\delta) = Min \{ \sigma_x(\delta): x \in D \}$, and we see
 that if $x \in \mathbb{H}^2$, $y, z \in S^1_{\infty}$,
 and $d_x(y,z) < \sigma(\delta)$, then
 $d_{\overline{\mathcal{I}} x}(\overline{\mathcal{I}} y,
 \overline{\mathcal{I}} z) < \delta$.
 It follows from the definition that
 $\sigma(\delta) \rightarrow 0$ as $\delta \rightarrow 0$.

\begin{lem} \label{distortion}
 Suppose $x \in \mathbb{H}^2$, $X \subset S^1_{\infty}$,
 $\epsilon >0$ such that
 $N_{\epsilon}(X, d_x) = S^1_{\infty}$,
 and $\delta \in \sigma^{-1} [\epsilon, 1]$.
 Then $N_{\delta}(\overline{\mathcal{I}} X, d_{\mathcal{I} x})
                             = S^2_{\infty}$.
\end{lem}

\begin{proof}
 Let $y \in S^2_{\infty}$, and let $z \in \overline{\mathcal{I}}^{-1} y$.
 Then $d_x(X, z) < \epsilon \leq  \sigma(\delta)$, and so
 $ d_{\mathcal{I} x}(\overline{\mathcal{I}}X, y)=
        d_{\mathcal{I} x}(\overline{\mathcal{I}}X, \overline{\mathcal{I}}z)
                  < \delta$.
\end{proof}

\subsection{Relating limit sets} ${}$

Now, let $g:S \rightarrow M$ be an immersion which is transverse
 to the suspension flow, and let $A_0(S) \subset A(S) \subset \mathbb{H}^2$
 be as given in Section \ref{immersions}.
 Let $G = g_* \pi_1 S$.
 Our next goal is to relate the limit set of $G$ to the limit set
 of $A(S)$.

 The map $\phi$ (defined in Section \ref{immersions})
 induces a homeomorphism between
 the universal cover of $S$ and the convex set
 $int \, A(S) \subset \mathbb{H}^2$.
 Thus, every element of $\pi_1 S$ acts on $int \,\, A(S)$, and
 $int \,\, A(S)/\pi_1 S \cong S$.  See \cite{CLR} for more details.
 For $\gamma \in \pi_1 S$, let $\tau_{\gamma}:A(S) \rightarrow A(S)$
 be the corresponding action.
 Note that, if $g_* \gamma = \mathcal{I}_* \gamma^{\prime}$,
 for some $\gamma^{\prime} \in \pi_1 F$, then
 $\tau_{\gamma} = \tau_{\gamma^{\prime}}|_{A(S)}$.

The set $A(S)$ may be decomposed as a union
 of connected sets $B_j$, where
 $B_j = \tau_{\delta_j} A_0(S)$, and $\delta_j \in \pi_1 S$.

\begin{lem} \label{limitpoints}
Fix $x_0 \in A_0(S)$.  For every point $y \in B_{j{\infty}}$, there
 is a sequence of elements $\gamma_i \in \pi_1 S$,
 such that $g_* \gamma_i \subset \mathcal{I}_* \pi_1 F$,
 and $\tau_{\gamma_i}  x_0 \rightarrow y$.
\end{lem}

\begin{proof}
 Let  $y \in  (B_j)_{\infty}$ for some $j$.
 Let $z = \overline{\tau_{\delta_j}}^{-1} y \in A_0(S)_{\infty}$.
 It follows from the definition of $A_0(S)$ that there is a sequence
 $\gamma_i^{\prime}$ in $\pi_1 S$, such that
 $\tau_{\gamma_i^{\prime}} x \rightarrow z$ for any $x \in A(S)$, with
 $g_*\gamma_i^{\prime} \in \overline{\mathcal{I}} \pi_1 F$ for all $i$.
 Let $\gamma_i = \delta_j \gamma_i^{\prime} \delta_j^{-1}$.
 Then
\begin{eqnarray*}
 g_* \gamma_i &=& (g_* \delta_j) \, (g_* \gamma_i^{\prime})\,
                                  (g_* \delta_j^{-1})\\
          &\in& (g_* \delta_j) \, \mathcal{I}_*\pi_1 F \, (g_*\delta_j^{-1})\\
         &=& \mathcal{I}_* \pi_1 F,
\end{eqnarray*}
and
\begin{eqnarray*}
\tau_{\gamma_i} x_0 &=& 
  \tau_{\delta_j} \tau_{\gamma_i^{\prime}} (\tau_{\delta_j}^{-1} x_0)\\
 &\rightarrow& \overline{\tau_{\delta_j}} z\\
                       &=&y
\end{eqnarray*}
\end{proof}

Let $G = g_*(\pi_1 S)$.  Then we have:

\begin{lem} \label{limitset}
$\overline{\mathcal{I}} A(S)_{\infty} \subset \Lambda(G)$.
\end{lem}

\begin{proof}
 Since $A(S)_{\infty} = \overline{\bigcup B_{j \infty}}$,
 it is enough to show that
 $\overline{\mathcal{I}}(y) \in \Lambda(G)$,
 for every $y \in  B_{j \infty}$.

 Fix $x_0 \in int \, A(S)$.  
 Then by Lemma \ref{limitpoints}, for every $y \in B_{j\infty}$,
 there is a sequence $\{\gamma_i\} \subset \pi_1 S$ such
 that $g_* \gamma_i \in \mathcal{I}_* \pi_1 F$ and 
 $\tau_{\gamma_i}(x_0) \rightarrow y$ in
 $\overline{\mathbb{H}}^2$.
 Let $z_0 = \mathcal{I} x_0 \in \mathbb{H}^3$,
 and then by the commutative diagram in Section \ref{hypgeom},
 we have 
\begin{eqnarray*}
g_*(\gamma_i) z_0 &=& \mathcal{I} \tau_{\gamma_i} \mathcal{I}^{-1} (z_0)\\
  &=& \mathcal{I} \tau_{\gamma_i} x_0\\
 &\rightarrow & \overline{\mathcal{I}}(y).
\end{eqnarray*}
Therefore, $\overline{\mathcal{I}} y \in \Lambda(G)$.
\end{proof}

\subsection{Thickness and sparsity} ${}$

 An observer inside $A(S)$ who is very far from $Frontier\, A(S)$
 will perceive $A(S)_{\infty}$ to occupy a very dense subset of his field
 of vision.  Our next goal is to quantify this statement.
  For every $x \in \mathbb{H}^2$ and $Y \subset S^1_{\infty}$, define
 the ``visual sparsity of $Y$ with respect to $x$'' to be:
 $$ \theta(Y,x) = \inf\{r: N_r(Y, d_x)= S^1_{\infty}\}.$$
 For a subset $X \subset \mathbb{H}^2$, we define
 the ``visual sparsity of $Y$ with respect to $X$'' to be:
 $$\theta(Y,X) = \sup \{ \theta(Y,x): x \in X \}.$$
 We have the following lemma:

 \begin{lem} \label{rd}
 For any subset $X \subset A(S)$, we have:
 $$d(Frontier A(S),X) = \rho(2\theta(A(S)_{\infty},X)),$$
 where $\rho$ is the function given in Definition \ref{rho}.
 \end{lem}

 Thus, if the distance from the frontier is large, then the visual sparsity
 of $A(S)_{\infty}$ is small.
 We shall need another definition in the proof.
 If $(Y, d)$ is a metric space, and $Z \subset Y$, then
 $Diam(Z,d)$ denotes the diameter of $Z$, with respect to the metric $d$.

\begin{proof}
 We shall assume that there is a component,
 $\beta$, of the frontier of $A(S)$ which is closest
 to $X$, and a point $x \in X$ which is closest to $\beta$;
 the proof in the general case follows from a continuity argument.

 Let $x \in X$ be the point closest to $\beta$.
 Form a triangle by adding geodesic
 rays from $x$ to the endpoints of $\beta$, and let $\theta_0$ be
 the angle of this triangle at $x$.
 Let $I_{\beta}$ be the interval in $S^1_{\infty} - A(S)_{\infty}$
 subtended by $\beta$, and note that $Diam(I_{\beta}, d_x) = \theta_0$.
 We have:
 \begin{eqnarray*}
\theta(A(S)_{\infty},X) &=& (1/2) \sup \{ Diam(I, d_y) : y \in X \textrm{ and }
 I \textrm{ is a component of }
        S^1_{\infty} - A(S)_{\infty} \}\\
 &=& (1/2) Diam(I_{\beta}, d_x)\\
 &=& (1/2) \theta_0,
\end{eqnarray*}
 and so
 $$d(Frontier A(S),X) = d(\beta,x) = \rho(\theta_0)
 = \rho(2\theta(A(S)_{\infty},X)).$$
 \end{proof}

\subsection{Thickness criterion} ${}$

 We are now in a position to prove the ``thickness criterion''.
Suppose that $g_i:S_i \rightarrow M$ is a sequence
 of quasi-Fuchsian, cut-and-cross-join immersions,
 as defined in Section \ref{immersions},
 and let $G_i = g_{i*} \pi_1 S_i$.
 Let $A_0(S_i) \subset A(S_i) \subset \mathbb{H}^2$ be as
 defined in Section \ref{immersions}.

\begin{lem} \label{thick}
 (Thickness criterion)
 If $d(A_0(S_i), Frontier A(S_i)) \rightarrow \infty$, then
 $t(G_i) \rightarrow \infty$.
\end{lem}

\begin{proof}
For each $i$, we have an associated finite cover $\widehat{F}_i$
 of $F$, with a pair of distinguished curves
 $\hat{\alpha}_i$ and $\widehat{f \alpha}_i$, and
 a distinguished regular neighborhood
 $N_i$ of $\widehat{\alpha}_i \cup \widehat{f \alpha}_i$.
  We let  $\widehat{F}_i^- = \widehat{F}_i - N_i$,
 which is identified with a subset $S_{i,1}$ of $S_i$,
 and we let $S_{i,2} = S_i - int \, S_{i,1}$.
 After a homotopy, we may assume that the maps $g_i|_{S_{i,2}}$ all agree;
 i.e. for any $i,j$ there is a homeomorphism
 $h_{ij}: S_{i,2} \rightarrow S_{j,2}$ such that
 $g_i|_{S_{i,2}} = g_j h_{ij}|_{S_{i,2}}$.
 Let $\widetilde{S}_i$ be the universal cover of $S_i$,
 which we may identify with $int \, A(S_i)$, and let
 $\tilde{g}_i: \widetilde{S}_i \rightarrow \mathbb{H}^3$
 be a lift of $g_i$.
 Let $\widetilde{S}_{i,1}$ be the full pre-image of $S_{i,1}$
 in $\widetilde{S}_i$.
 Since the maps $g_i|_{S_{i,2}}$ all agree, it follows that
 there is a number $R$ such that
 $N_R(\tilde{g}_i(\widetilde{S}_{i,1})) \supset \tilde{g}_i \widetilde{S}_i$
 for every $i$. 
 
 Let $D_i$ be the set of points in $A_0(S_i)$ which project to
 $S_{i,1}$. 
 Note that $G_i(\tilde{g}_i D_i) \supset \tilde{g}_i \widetilde{S}_{i,1}$
 and so $N_R(G_i(\tilde{g}_i D_i)) \supset \tilde{g}_i \widetilde{S}_i$.
 Therefore, by Lemma \ref{technical}, it is enough to prove the following:\\
\\
 \textit{Claim:}  There is a sequence
 $\delta_i \rightarrow 0$ such that
 $N_{\delta_i}(\Lambda_i, d_{\tilde{g}_i x}) = S^2_{\infty}$
 for all $x \in D_i$.\\
\\
\textit{Proof of Claim:}
 Let $x \in A_0(S_i)$, let $r_i = d(A_0(S_i), Frontier \, A(S_i))$,
 and let $\theta_i = \theta(A(S_i)_{\infty},A_0(S_i))$.
 By the definition of $\theta_i$, we have
 $N_{\theta_i}(A(S_i)_{\infty}, d_x) = S^1_{\infty}$.
 By Lemma \ref{rd}, $\rho(2 \theta_i) = r_i$.
 Since $r_i \rightarrow \infty$, then, by the definition
 of $\rho$, we see that $\theta_i \rightarrow 0$. 
  Recall the function $\sigma$ defined earlier in this section,
 and define $\delta_i = \inf \{ \sigma^{-1} [\theta_i, 1] \}$.
 Since $\sigma(x) \rightarrow 0$ as $x \rightarrow 0$,
 then $\delta_i \rightarrow 0$. 

 We have:
\begin{eqnarray*}
 N_{\delta_i}(\overline{\mathcal{I}} A(S_i)_{\infty}, d_{\tilde{g}_i x})
 &=& N_{\delta_i}(\overline{\mathcal{I}} A(S_i)_{\infty}, d_{\mathcal{I}x})
 \,\,\,\,\,\textrm{(since } x \in A_0(S_i) \textrm{)},\\
 &=& S^2_{\infty},  \textrm{ by Lemma \ref{distortion}}.
\end{eqnarray*}

 By Lemma \ref{limitset},
 $\Lambda_i \supset \overline{\mathcal{I}} A(S_i)_{\infty}$,
 and so $N_{\delta_i}(\Lambda_i, d_{\tilde{g}_i x}) = S^2_{\infty}$.
 This proves the claim and hence the lemma.
 \end{proof}

\section{Proof of Theorem \ref{main}}

\subsection{Outline} ${}$

 First we explain the idea for constructing
 cross-join surfaces which are quasi-Fuchsian.
  Replacing the given monodromy $f:F \ra F$ with a power,
 we may assume that it has a fixed point.  Then we fix
 a lift ${\wt f}: \mathbb{H}^2 \ra \mathbb{H}^2$,
 with a fixed point, $p$.
   We choose a stable leaf $R \supset \{p \}$
 of the invariant singular foliation of $\mathbb{H}^2$.
 Choose a fixed simple closed curve $\a$ in $F$,
 and consider a lift  ${\wt \a}$
 of $\a$ which intersects $R$.
 We wish to construct a cross-join surface
 $S$ so that the curves ${\wt \a}$ and  ${\wt f} {\wt \a}$
 are boundary components of $A_0(S)$.
 For then the map ${\wt f} = \tau_{0,j}$ for some $j$,
 and since ${\wt f}$ has a fixed point, the surface
 $S$ will be quasi-Fuchsian by Lemma \ref{coarse}.

 The construction of such a surface $S$ reduces to
 an immersion-to-embedding problem in the fiber $F$.
 Indeed, we require the existence of\\
\\
I. a finite cover ${\wh F} \ra F$, containing embedded, disjoint lifts
 ${\wh \a}$ and ${\wh f \a}$, which lift to ${\wt \a}$
 and ${\wt f} { \wt \a}$ under the induced covering
 of ${\wh F}$ by $\mathbb{H}^2$.\\
\\
 To find a \textit{thick}, quasi-Fuchsian cross-join surface $S$,
 it is enough to show that, in addition, $A_0(S)$ is far from
 the boundary of $A(S)$ (by Lemma \ref{thick}).
 For this, it is enough to show that
 the curves ${\wt f}^{-1} {\wt \a}, {\wt f} \a, {\wt f} {\wt \a},
 {\wt f}^2 \a$ project to curves which are all far apart in ${\wh F}$.
 So we must first ensure that\\
\\
II. the curves ${\wt f}^{-1} {\wt \a}, {\wt f} \a, {\wt f} {\wt \a},
 {\wt f}^2 \a$ are far apart in $\mathbb{H}^2$,\\
\\
 and then show that\\
\\
III. the projections of  ${\wt f}^{-1} {\wt \a}, {\wt f} \a, {\wt f} {\wt \a},
 {\wt f}^2 \a$ are far apart in ${\wh F}$,
 and are homeomorphic lifts of $f^{-1} \a, ..., f^2 \a$.\\

 For Condition II, it is enough to understand the dynamics
 of the pseudo-Anosov map ${\wt f}$, when restricted to the stable leaf $R$.
 This is a straight-forward matter. The analysis to this point is all carried
 out in $\mathbb{H}^2$, and is contained in Sub-section \ref{h2}.

 Condition III (which includes Condition I) reduces
 to an immersion-to-embedding problem, whose statement requires
 some more notation.
 Corresponding to each curve ${\wt \a}_i$, there is a group element
 $\g_i \in \pi_1 F$, and it is convenient to replace the
 curves ${\wt f}^i {\wt \a}$ with geodesics $Ax(\gamma_i)$. 
 If ${\wh F}$ satisfies Condition III,
 then there is an intermediate (non-compact) surface
 $F^* = \mathbb{H}^2/<\gamma_{-1}, ..., \gamma_2>$
 covering ${\wh F}$. We must first ensure that  
  the projection of $Ax(\g_{-1}) \cup ... \cup Ax(\g_2)$ 
 has a large, embedded collar in the intermediate cover $F^*$
 (Lemma \ref{pingpong}).  This is the content of Sub-section \ref{middle}.

 The collar immerses from $F^*$ into $F$ under a covering map.
 To satisfy Conditions I and III (and complete
 the proof of Theorem \ref{main}),
 it is enough to lift this immersion
 to an embedding in a finite cover of $F$.
 This is done in Lemma \ref{dist},
 and is the main content of Sub-section \ref{end}.

\subsection{Arranging things in $\mathbb{H}^2$} \label{h2} ${}$
 
 We now return to the situation of Theorem \ref{main}.
 After passing to a finite cover, we 
 have a fixed hyperbolic 3-manifold $M$,
 fibering over $S^1$, with pseudo-Anosov monodromy $f:F \rightarrow F$.
 Since $f$ is pseudo-Anosov, it has periodic points in $F$, and thus,
 after passing to a further cyclic cover of $M$, we may replace
 $f$ with a power, and assume that $f$ has a fixed point $p$. We fix
 a lift $\widetilde{f}$ of $f$ to $\mathbb{H}^2$ with a fixed point
 $\widetilde{p}$,
 let $\pi:\mathbb{H}^2 \rightarrow F$ be the universal covering map, and
 let $R$ be a geodesic ray 
 with an endpoint at $\widetilde{p}$, which is an ``attracting axis''
 for $\tilde{f}$-- i.e. $\tilde{f} R = R$, and $\tilde{f}$ shrinks distances
 on $R$ by a factor of $1/\lambda$, where $\lambda$ is the dilatation of $f$.

\begin{lem} \label{distance}
 Let $p_i$ be a sequence of points on $R$ such that
 $d(p_i, \widetilde{p}) \rightarrow \infty$, and let $j, k$ be distinct
 integers.
 Then $d(\widetilde{f}^j p_i, \widetilde{f}^k p_i) \rightarrow \infty$.
\end{lem}

\begin{proof}
 The map $\widetilde{f}$ acts on $R$
 by the formula
 $d(\widetilde{f} x, \widetilde{p}) = \lambda^{-1} d(x, \widetilde{p})$.
 The lemma follows.
\end{proof}

\begin{lem} \label{angle}
Let $\alpha$ be any simple closed geodesic in $F$.
 Then the angles of intersection between $\pi R$ and $\alpha$ are bounded
 below by some angle $\theta > 0$.
\end{lem}

\begin{proof}
 Note that there is a well-defined intersection angle at each
 point of $\alpha \cap \pi R$, since $\alpha$ and $\pi R$ are simple
 geodesics.
 Suppose there is no lower bound on these angles.
 Then there is an infinite sequence of
 points $q_i \in \pi R \cap \alpha$,
 with angles of intersection tending to 0.
 Let $q$ be an accumulation point of this sequence.

 Let $\Lambda \subset F$ be the invariant
 geodesic lamination for $f$. Then the accumulation points of
 $\pi R$ are all contained in $\Lambda$, so
 $\Lambda \cup \pi R$ is a closed subset of $F$.
 Therefore $q \in \alpha \cap (\Lambda \cup \pi R)$,
 and by continuity, the angle of intersection between
 $\alpha$ and $(\Lambda \cup \pi R)$ must be zero at $q$.
 Therefore $\alpha$ coincides with either a leaf of $\Lambda$
 or else with $\pi R$.  However these are both impossible,
 since $\pi R$ is not closed, and $\Lambda$ contains no closed leaves.
\end{proof}

\begin{lem} \label{int}
Every primitive class in $H_1(F, \mathbb{Z})$ has a geodesic representative
 whose intersection with $\pi R$ is infinite.
\end{lem}

\begin{proof}
For every primitive class there is a pair of simple, closed, geodesic
 representatives,
 $\alpha_1, \alpha_2$, intersecting transversely, so that the complement
 of their union is a collection of open disks, whose
 diameter is bounded by some constant $D$.  Since $\pi R$ is
 a geodesic ray, every sub-segment of $\pi R$ of length greater than $D$ must
 intersect one of the two representatives.  Since $\pi R$ is an infinite
 ray, then we conclude that
 $\pi R \cap (\alpha_1 \cup \alpha_2)$ is infinite.
 Therefore either $\pi R \cap \alpha_1$
 or $\pi R \cap \alpha_2$ is infinite.
 \end{proof}

By Lemma \ref{int}, we may assume that $\alpha \cap \pi R$ is infinite.

\begin{lem} \label{points}
There is an infinite sequence of points $p_i \subset R$, such that
 $\pi p_i \in \alpha$ for all $i$, and
 $d(\tilde{p}, p_i) \rightarrow \infty$.
\end{lem}

\begin{proof}
 By assumption, there is an infinite sequence of points $\{p_i\} \subset R$
 which project to $\alpha$.
 Since the intersection of any finite sub-ray of the geodesic $\pi R$
 with the geodesic $\alpha$ is finite,
 we may pass to a subsequence of $\{p_i\}$, and
 assume that $d(p_i, \tilde{p}) \rightarrow \infty$.
\end{proof}

 The free homotopy class of $f^j \alpha$ corresponds to a certain
 conjugacy class in $\pi_1 F$;
 we shall now specify a representative in $\pi_1 F$ for each such class.
  Let $\widetilde{\alpha}$ be a fixed lift of
 $\alpha$ to $\mathbb{H}^2$ which intersects $R$, and 
 fix the basepoint for $\pi_1 F$ at $p$. There is a unique oriented arc $A_j^+$
 from $\widetilde{p}$ to $\widetilde{f}^j \widetilde{\alpha}$ along the geodesic ray $R$,
 and $\pi A_j$ is an arc from $p$ to $f^j \alpha$; let
 $\gamma_j \in \Gamma$ be
 represented by the loop $A_j^+ * f^j \alpha * A_j^-$.

 If $\gamma \in Isom^+(\mathbb{H}^2)$ is a hyperbolic isometry,
 let $Ax(\gamma)$ denote the axis of $\gamma$. 
 Let $s_0 = d(\widetilde{\alpha}, \tilde{p})$.

\begin{lem} \label{distance2}
Given  any number $J > 0$, there is a number $L > 0$,
 such that, if $s_0 > L$, then $d(Ax(\gamma_j), Ax(\gamma_k))>J$
 for all distinct $j,k \in \{ -1, .., 2\}$.
\end{lem}

\begin{proof}
 We first define some useful constants.
 By construction, for every integer $i$, there is a number $K_i$ such that
 $\widetilde{f}^i \widetilde{\alpha}$ is contained
 in the $K_i$-neighborhood of $Ax(\gamma_i)$.
 Let $K = Max\{ K_{-1}, ..., K_2 \}$.
 Also, since $\pi Ax(\gamma_i)$ is freely homotopic to the simple
 closed curve $f^i \alpha$, then $\pi Ax(\gamma_i)$
 is itself simple, and therefore,
 by Lemma \ref{angle}, the  angle of intersection
 made by $R$ with $Ax(\gamma_i)$ is greater than
 some constant $\theta_i$, independent of the original choice
 of $\widetilde{\alpha}$, and we let
 $\theta = Min \{ \theta_{-1}, ..., \theta_2 \}$.
 
 Let $q_i = Ax(\gamma_i) \cap R$.  Note that
 $q_0 = Ax(\gamma_0) \cap R = \tilde{\alpha} \cap R$.
 Since $N_K(Ax(\gamma_i)) \supset \widetilde{f}^i \widetilde{\alpha}$,
 then $d(\widetilde{f}^i q_0, Ax(\gamma_i))
             \leq d(\widetilde{f}^i \tilde{\alpha}, Ax(\gamma_i)) \leq K$.
 Therefore, by the hyperbolic law of sines, we have
\begin{eqnarray*}
 d(\widetilde{f}^i q_0, q_i) \leq K^{\prime},
\end{eqnarray*}
 where $K^{\prime} = arcsinh(sinh(K)/sin(\theta))$.
 So, for distinct integers $j,k \in \{-1, .., 2 \}$ we have:
\begin{eqnarray*}
d(q_j, q_k) &\geq& d(\widetilde{f}^j q_0, \widetilde{f}^k q_0)
     - d(\widetilde{f}^j q_0, q_j) - d(\widetilde{f}^k q_0, q_k) \\
&\geq& d(\widetilde{f}^j q_0, \widetilde{f}^k q_0)
 - 2 K^{\prime}\\
&\rightarrow& \infty,\textrm{ as } s_0 \rightarrow \infty,
 \textrm{ by Lemma }\ref{distance}.
\end{eqnarray*}

 Since
 $d(q_j, q_k) \rightarrow \infty$, and since the
 angles of intersection of $R$ with $Ax(\gamma_j)$
 and $Ax(\gamma_k)$ are bounded below by $\theta$,
 we conclude that
 $d(Ax(\gamma_j), Ax(\gamma_k)) \rightarrow \infty$,
 as $s_0 \rightarrow \infty$.
 \end{proof}

\subsection{Arranging things in the intermediate cover} \label{middle} ${}$

 Let $F^* = \mathbb{H}^2/<\gamma_{-1}, \gamma_0, \gamma_1, \gamma_2>$,
 and let $\nu: \mathbb{H}^2 \rightarrow F^*$ be the universal cover.
 Let $\alpha_j^* = \nu Ax(\gamma_j)$. 

If  $\mathcal{B} = \beta_1 \cup ... \cup \beta_n$
 is the union of a collection of disjoint, simple, closed
 geodesics in a surface $S$, 
 we define the \textit{tube radius} of $\mathcal{B}$
 to be the supremum of all numbers $r$ such that
 $\overline{N_r}(\mathcal{B})$ consists of $n$ disjoint, properly
 embedded annuli in $S$.  We denote this quantity by $tuberad(\mathcal{B})$.

 \begin{lem} \label{pingpong}
 Given any number $J > 0$, there is a number $L>0$ such that,
 if $s_0 > L$, then $\alpha_{-1}^*, ..., \alpha_2^*$ is a collection
 of disjoint, simple closed geodesics, and
 $tuberad(\alpha_{-1}^* \cup ... \cup \alpha_2^*) > J$.
\end{lem}

\begin{proof}
 Let $Arc_i \subset Ax(\gamma_i)$ be a fundamental domain for the action
 of $\gamma_i$ on $Ax(\gamma_i)$.
 It suffices to prove that, if $s_0$ is large enough, then for any
 $i,j \in \{-1, ..., 2 \}$ and
 any $\gamma \in <\gamma_{-1}, ..., \gamma_2>$,
 either $d(Arc_i,\gamma Ax(\gamma_j)) > J$,
 or else $i=j$ and $\gamma \in <\gamma_i>$.

 Then let $R_i^{\pm}$ be a pair of geodesics orthogonal to
 $Ax(\gamma_i)$, equidistant from $Ax(\gamma_i) \cap R$,
 such that $\gamma_i R_i^- = R_i^+$.
 Let $U_i^{\pm}$  be the component of
 $S^1_{\infty} - R_i^{\pm}$ containing
 the attracting fixed point of $\gamma_i^{\pm 1}$. 
 See Figure \ref{ping}.

\begin{figure}[!ht] \label{ping}
{\epsfxsize=4in \centerline{\epsfbox{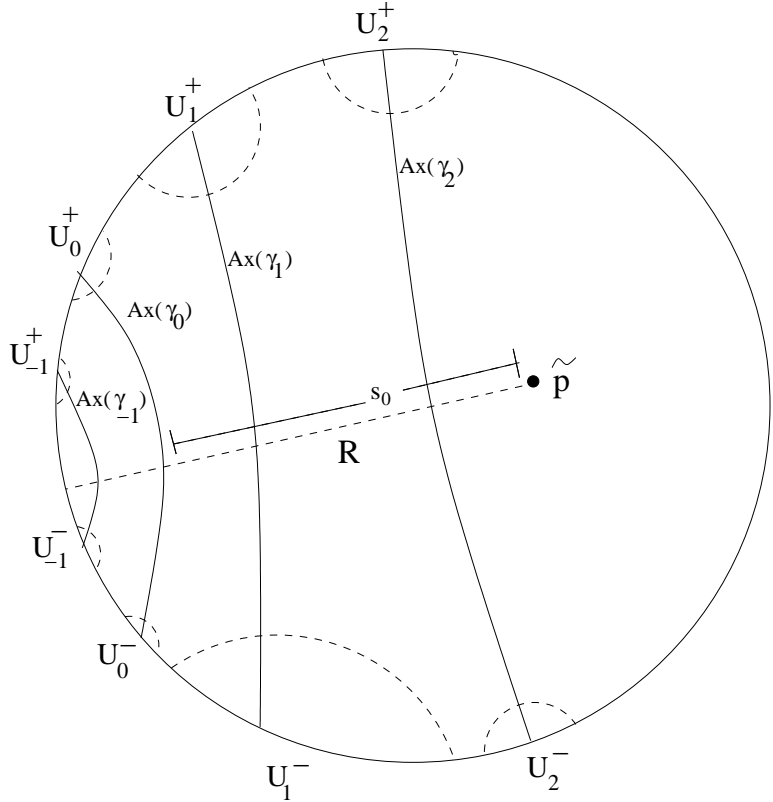}}\hspace{10mm}}
\caption{Increasing the tube radius.}
\end{figure}
 
 By increasing $s_0$ and applying Lemma \ref{distance2},
 we may assume that $d(Ax(\gamma_i), Ax(\gamma_j))$
 is large enough so that
 $(U_i^+ \cup U_i^-) \cap (U_j^+ \cup U_j^-) = \emptyset$
 for all distinct $i,j \in \{-1, ..., 2\}$. 
 Then we have: 

 \begin{eqnarray*}
&& \gamma_i^{\pm 1} U_i^{\pm} \subset U_i^{\pm}, \textrm{ and}\\
&& \gamma_j^{\pm 1} (U_i^+ \cup U_i^-) \subset U_j^{\pm},
 \textrm{ for any distinct }i,j \in \{-1, ..., 2\}.
\end{eqnarray*}
 Similarly, by increasing $s_0$, we may also assume that:\\

 \textit{a.}  for any $i,j \in \{-1, ..., 2 \}$, any $k \neq j$, 
 and any non-zero integer $n$,
 if $\beta$ is a geodesic with endpoints in $\gamma_k^n(U_j^+ \cup U_j^-)$,
 then $d(Arc_i, \beta) > J$.

 \textit{b.} $d(Ax(\gamma_i), Ax(\gamma_j)) > J$
 for all distinct $i,j \in \{-1, ..., 2\}$.\\

 Let $\gamma
 \in <\gamma_{-1}, ..., \gamma_2>$, and let $i,j \in \{-1, ..., 2\}$.
 Putting $\gamma$ in normal form, we have
 $\gamma = \gamma_{s_1}^{t_1} ... \gamma_{s_k}^{t_k}$, where each consecutive
 pair of subscripts is distinct. We must consider several cases.\\
\\
 If $k=1$, and $s_1 \neq j$:\\
 Then $\gamma = \gamma_{s_1}^{t_1}$, and since
 $\gamma Ax(\gamma_j)= \gamma_{s_1}^{t_1} Ax(\gamma_j)$ has endpoints in
 $\gamma_{s_1}^{t_1}(U_j^+ \cup U_j^-)$,
 we have $d(Arc_i, \gamma Ax(\gamma_j)) > J$, by assumption \textit{a}.\\
\\
 If $k=1$, $s_1 = j$ and $i \neq j$:\\
 Then $d(Arc_i, \gamma Ax(\gamma_j)) = d(Arc_i, Ax(\gamma_j)) >J$,
 by assumption \textit{b}.\\
\\
 If $k>1$:\\
 The endpoints of $\gamma Ax(\gamma_j)$ are
 both contained in $\gamma_{s_1}^{t_1}(U_{s_2}^+ \cup U_{s_2}^-)$.
 Since $s_1 \neq s_2$,  we may apply assumption \textit{a},
 to conclude that  $d(Arc_i, \gamma Ax(\gamma_j)) > J$.\\
\\
Thus we have shown that $d(Arc_i, \gamma Ax(\gamma_j)) > J$,
 unless $k=1$, and $s_1 = i = j$.
 Equivalently,  $d(Arc_i, \gamma Ax(\gamma_j)) > J$,
 unless  $i = j$, and $\gamma \in <\gamma_i>$, which is what we had to show.
 \end{proof}

\subsection{Arranging things in the finite cover} \label{end} ${}$

\begin{lem} \label{dist}
The cover $F^* \rightarrow F$ factors through
 a finite cover $\widehat{F} \rightarrow F$,
 as in the following diagram:
\begin{eqnarray*}
F^*   & \psi \atop \rightarrow & \widehat{F}\\
      &    \searrow &   \downarrow \\  
      & & F,
\end{eqnarray*}
 such that\\
i.  $\psi \alpha_{-1}^*, ..., \psi \alpha_2^*$ are embedded
  and pair-wise disjoint, and\\
ii. $tuberad(\alpha_{-1}^* \cup ... \cup \alpha_2^*)
       = tuberad(\psi \alpha_{-1}^* \cup ... \cup \psi \alpha_2^*)$.
\end{lem}

\begin{proof}
 Let $r = tuberad(\alpha_{-1}^* \cup ... \cup \alpha_2^*)$.
 Using notation from the proof of Lemma \ref{pingpong}, let
 $$X_{j,k} =  \{ \gamma \in \pi_1 F - \pi_1 F^*|
 N_r(Arc_k) \cap \gamma N_r(Ax(\gamma_j)) \neq \emptyset \}.$$
 By discreteness, $X_{j,k}$ is finite.
 Let $X = \bigcup_{j,k \in \{-1, ..., 2\}} X_{j,k}$.
 Since $\pi_1 F$ is LERF (by \cite{S}), there is a finite-index subgroup
 of $\pi_1 F$ which contains $\pi_1 F^*$, but is disjoint from $X$.
 The corresponding finite cover $\widehat{F} \rightarrow F$
 satisfies the hypotheses of the lemma.
\end{proof}

 Now let $s_i$ be a sequence of numbers approaching $\infty$.
 By Lemmas \ref{points} and \ref{pingpong}, for every $i$, there is a choice
 of $\widetilde{\alpha}$ such that
 $tuberad(\alpha_{-1}^* \cup ... \cup \alpha_2^*) > s_i$.
 By Lemma \ref{dist}, the corresponding
 cover $F^*_i$ factors through a finite cover $\widehat{F}_i \rightarrow F$,
 with induced covering map $\psi: F^*_i \rightarrow \widehat{F}_i$,
 such that $tuberad(\psi \alpha_{-1}^* \cup ... \cup \psi \alpha_2^*) > s_i$.

  Let $\widehat{f^j \alpha} \in \widehat{F}_i$ be the unique lift
 of $f^j \alpha$ in the free homotopy class of 
 the geodesic $\psi \alpha_j^*$.
 Then, using the notation of Section \ref{immersions},
 we let $S_i = S(\widehat{\alpha}, \widehat{f \alpha})$,
 and let  $g_i:S_i \rightarrow M$ be the corresponding immersion.
 
\begin{lem}
The immersion $g_i:S_i \rightarrow M$ is quasi-Fuchsian.
\end{lem}

\begin{proof}
In the notation of Section \ref{immersions}, consider the subset
 $A(S_i) = \bigcup_j A_j$ of $\mathbb{H}^2$.
 Note that $\widetilde{\alpha}$ and $\widetilde{f} \widetilde{\alpha}$ are
 both components of $Frontier \, A_0(S_i)$.  Therefore we may suppose
 $\widetilde{f} = \tau_{0j}$ for some $j$. 
 By construction, $\widetilde{f}$ has a fixed point, $\widetilde{p}$.
 Therefore, by Lemma \ref{coarse}, the immersion
 $g_i:S_i \rightarrow M$ is quasi-Fuchsian.
\end{proof}

As in the proof of Lemma \ref{distance2}, we
 may choose a number $K$, independent of $i$, such that,
 for any $j \in \{-1, ...,2\}$, the curve $\widehat{f^j \alpha}$
 is contained in the $K$-neighborhood of $\psi \alpha_j^*$.

\begin{lem} \label{dist2}
$d(A_0(S_i), Frontier \, A(S_i)) > s_i -2K$.
\end{lem}

\begin{proof}
 Suppose there is a geodesic path $\beta$ from $A_0(S_i)$
 to $Frontier \, A(S_i)$, of length at most $s_i-2K$.
 Note that some sub-arc of $\beta$ projects to a
 curve connecting either $\widehat{f \alpha} \cup \widehat{f^2 \alpha}$,
 or $\widehat{\alpha} \cup \widehat{f^{-1} \alpha}$.
 So either $d(\widehat{f \alpha}, \widehat{f^2 \alpha}) \leq s_i -2K$
 or  $d(\widehat{\alpha},\widehat{f^{-1} \alpha}) \leq s_i -2K$.
 Since $\widehat{f^j \alpha}$
 is contained in the $K$-neighborhood of $\psi \alpha_j^*$ for any
 $j \in \{-1, ..., 2\}$, we conclude that
 $tuberad(\psi \alpha_{-1}^* \cup ... \cup \psi \alpha_2^*) \leq s_i$,
 contradicting the choice of $s_i$.
\end{proof}

 Let $G_i = Image(g_{i*})$.

\begin{lem} \label{done}
$t(G_i) \rightarrow \infty$.
\end{lem}

\begin{proof}
This is a consequence of Lemmas \ref{thick} and \ref{dist2}.
\end{proof}

The proof of Theorem \ref{main} is now complete.

\vspace{5mm} \noindent Mathematics Department\\
 SUNY at Buffalo\\
Buffalo, NY 14260-2900\\
 jdmaster@buffalo.edu


\begin{thebibliography}{999999}
\bibitem{CT} J. Cannon and W. P. Thurston, ``Group-invariant Peano
 Curves'', pre-print.

\bibitem{CL} D. Cooper and D. D. Long, ``Some surface subgroups survive
 surgery'', textit{Geom. Topol.} \textbf{5} (2001), 347-367.

\bibitem{CLR} D. Cooper, D. D. Long
 and A. W. Reid, ``Bundles and finite foliations'',
 \textit{Invent. Math.} \textbf{118} (1994), 255-283.

 
 
\bibitem{EM} D. B. A. Epstein, A. Marden,
``Convex hulls in hyperbolic space, a theorem of Sullivan,
 and measured pleated surfaces'',
 \textit{Analytical and geometric aspects of hyperbolic space
 (Coventry/Durham 1984)}, 113--253,
 London Math. Soc. Lecture Note Ser., 111,
 Cambridge Univ. Press, Cambridge, 1987.  


\bibitem{F} S. Fenley,
 ``Surfaces transverse to pseudo-Anosov flows and virtual
 fibers in $3$-manifolds'', \textit{Topology} \textbf{38} (1999),
 no. 4, 823--859. 

\bibitem{M} B. Mangum,
 ``Incompressible surfaces and pseudo-Anosov flows'',
 \textit{Topology Appl.} \textbf{87} (1998), no. 1, 29-51.

\bibitem{S} G. P. Scott,
 ``Subgroups of surface groups are almost geometric'',
\textit{J. London Math. Soc.} (2) \textbf{17} (1978), no. 3, 555--565. 

\bibitem{Wu} Ying-Qing Wu, ``Immersed surfaces and Dehn surgery,''
\textit{Topology} \textbf{43} (2004), 319-343. 


\end{thebibliography}
\end{document}